\documentstyle[12pt,amssymb]{article}
\begin{document}

{\Large

\noindent{\bf On factorization of $q$-difference equation }

\vspace{0.3cm}

\noindent{\bf for continuous $q$-Hermite polynomials} }

\bigskip

\noindent{\bf M. N. Atakishiyev}

\noindent {Instituto de Matem\'aticas, Universidad Nacional Aut\'onoma
de M\'exico \\ CP 62210 Cuernavaca, Morelos, M\'exico}

\noindent E-mail: mamed@matcuer.unam.mx

\bigskip

\noindent{\bf A. U. Klimyk}

\noindent Bogolyubov Institute for Theoretical Physics, Kiev 03143, Ukraine

\noindent E-mail: aklimyk@bitp.kiev.ua

\bigskip

\begin{abstract}

We argue that a customary $q$-difference equation for the continuous
$q$-Hermite polynomials $H_n(x|\,q)$ can be written in the factorized
form as $\left[\left({\cal D}_x^{\,q}\right)^2 -1 \right]\, H_n(x|\,q)\\
= \left(q^{-n}-1\right)\,H_n(x|\,q)$, where ${\cal D}_x^{\,q}$ is some
explicitly known $q$-difference operator. This means that the polynomials
$H_n(x|\,q)$ are in fact governed by the $q$-difference equation
${\cal D}_x^{\,q} \,H_n(x|\,q)= q^{-n/2}\,H_n(x|\,q)$, which is simpler
than the conventional one.

\end{abstract}

\medskip


\bigskip

It is well-known that the theory of $q$-Hermite polynomials is one
of the main instruments for studying $q$-oscillators and their
applications (see, for example, \cite{Mac, Bied}). In particular,
$q$-difference equations for $q$-Hermite polynomials are known to
be intimately connected with Hamiltonians of the corresponding
systems.

The  continuous $q$-Hermite polynomials are also related to the
oscillator representations (a special case of the discrete series
representations) of the quantum algebra ${\rm su}_q(1,1)$, which
are constructed with the aid of the creation and annihilation
operators for the $q$-oscillator (see, for example, \cite{KD}).

The aim of our paper is to study $q$-difference equations for the
continuous $q$-Hermite polynomials for both cases when $0<q<1$ and
$q>1$. The  continuous $q$-Hermite polynomials for $q>1$ are essentially
different from those for $0<q<1$.

The continuous $q$-Hermite polynomials of Rogers, $H_n(x|q)$, $0<q<1$,
are orthogonal on the finite interval $-1\leq x := \cos \theta \leq 1$,
$$
\frac{1}{2\pi}\,\int_{-1}^{1}\,H_m(x|\,q)\,H_n(x|\,q){\widetilde w}
(x|\,q)\,dx = \frac{\delta_{mn}}{(q^{n+1};q)_{\infty}}\,, \eqno(1)
$$
with respect to the weight function (we employ standard notations of the
theory of special functions, see, for example, \cite{GR} or \cite{AAR})
$$
\widetilde w(x|\,q):= \,\frac{1}{\sin \theta}\,
\left(\,e^{2{\rm i}\theta},\,e^{-2{\rm i}\theta}; q \right)_\infty \, . \eqno(2)
$$
These polynomials satisfy the $q$-difference equation
$$
D_{q}\,\left[{\widetilde w}(x|\,q)\,D_q\,H_n(x|\,q)\right] =
\frac{4\,q\,(1-q^{-n})}{(1-q)^2}\,H_n(x|\,q)\,{\widetilde w}(x|\,q)\,, \eqno(3)
$$
written in a self-adjoint form \cite{NSU}. The $D_q$ in (3) is the conventional
notation for the Askey-Wilson divided-difference operator defined as
$$
D_{q}\, f(x) := \frac{\delta_{q}\,f(x)}{\delta_{q}\,x}\,,  \eqno(4)
$$
$$
\delta_q\,g(e^{{\rm i}\,\theta}):= g(q^{1/2}\,e^{{\rm i}\,\theta}) -
g(q^{-1/2}\,e^{{\rm i}\,\theta})\,, \quad f(x)\equiv g(e^{{\rm i}\theta})\,,
\quad x=\cos\theta \,.
$$
In what follows we find it more convenient to employ the explicit expression
$$
D_{q}\,f(x)= \frac{\sqrt{q}}{{\rm i}(1-q)}\,\frac{1}{\sin\theta}
\left( e^{\,{\rm i}\ln q^{1/2}\,\partial_{\theta}}- e^{-{\rm i}\ln q^{1/2}\,
\partial_{\theta}}\right)f(x),\quad \quad \partial_{\theta}
\equiv\frac{d}{d\theta}\,,                   \eqno(5)
$$
for the $D_q$ in terms of the shift operators (or the operators of the
finite displacement, \cite{LL}) $e^{\pm a\,\partial_{\theta}}\,g(\theta)
:= g(\theta \pm a)$ with respect to the variable $\theta$. Although it
is customary to represent $q$-difference equation for the $q$-Hermite
polynomials in the self-adjoint form (3) (see \cite{KS}, p. 115), one
may eliminate the weight function $\widetilde w(x|\,q)$ from (3) by
utilizing its property that
$$
\exp \left(\pm \,{\rm i}\ln q^{1/2}\,\partial_{\theta}\right){\widetilde w}(x|\,q)
= - \frac{e^{\pm 2{\rm i}\theta}}{\sqrt{q}}\,{\widetilde w}(x|\,q)\,.    \eqno (6)
$$
The validity of (6) is straightforward to verify upon using the explicit
expression (2) for the weight function ${\widetilde w}(x|\,q)$.

Thus, combining (3) and (6) results in the $q$-difference equation
$$
\frac{1}{2{\rm i}\sin\theta}\,\left[\,\frac{e^{{\rm i}\theta}}
{1-q\,e^{-2{\rm i}\theta}}\,\left(e^{{\rm i}\ln q\,\partial_{\theta}}\,
- 1\right)+ \,\frac{e^{-{\rm i}\theta}}{1 - q\,e^{2{\rm i}\theta}}\,
\left(1 - e^{-{\rm i}\ln q\,\partial_{\theta}}\,\right)\right]\,H_n(x|\,q)$$
$$ = \left( q^{-n}-1\right)\,H_n(x|\,q)\,                           \eqno(7)
$$
for the continuous $q$-Hermite polynomials $H_n(x|\,q)$, which does not
explicitly contain the weight function ${\widetilde w}(x|\,q)$.

In connection with equation (7) it should be remarked that Koornwinder \cite{K} has
recently studied in detail raising and lowering relations for the Askey-Wilson
polynomials $p_n(x;a,b,c,d|\,q)$. We
recall that the Askey-Wilson fami\-ly for $a=b=c=d=0$ is known to reduce to
the continuous $q$-Hermite polynomials $H_n(x|\,q)$. So, as a consistency check,
one may verify that (7) is in complete agreement with particular case of the
equation $D\,p_n=\lambda_n\,p_n$ (i.e., equation (4.5) in \cite{K}) for
Askey-Wilson polynomials  with vanishing parameters $a,b,c,d$.

We are now in a position to show that equation (7) admits a factorization. Indeed,
with the help of two simple trigonometric identities
$$
\frac{e^{\,\pm\,{\rm i}\theta}}{2{\rm i}\sin\theta} =
\pm \,\frac{1}{1-e^{\,\mp\, 2{\rm i}\theta}}
$$
one can represent the left side of (7) as
$$
\frac{1}{2{\rm i}\sin\theta}\,\left(\,\frac{e^{{\rm i}\theta}}
{1-q\,e^{-2{\rm i}\theta}}\,\,e^{{\rm i}\ln q\,\partial_{\theta}}\,
- \,\frac{e^{-{\rm i}\theta}}{1 - q\,e^{2{\rm i}\theta}}\,
\,e^{-{\rm i}\ln q\,\partial_{\theta}}\, - \frac{e^{{\rm i}\theta}}
{1-q\,e^{-2{\rm i}\theta}}+ \,\frac{e^{-{\rm i}\theta}}{1 -
q\,e^{2{\rm i}\theta}}\right) H_n(x|\,q)$$
$$= \left[\, \frac{1}{1- e^{-2{\rm i}\theta}}\,\,e^{{\rm i}\ln q^{1/2}
\partial_{\theta}}\, \frac{1}{1- e^{-2{\rm i}\theta}}\,\, e^{{\rm i}\ln q^{1/2}
\partial_{\theta}} + \frac{1}{1- e^{2{\rm i}\theta}}\,\,e^{-{\rm i}\ln q^{1/2}
\partial_{\theta}}\, \frac{1}{1- e^{2{\rm i}\theta}}\,\,e^{-{\rm i}\ln q^{1/2}
\partial_{\theta}}\right. $$
$$
+ \left. \frac{q(1+q)}{(1+q)^2-4qx^2} - 1 \right] H_n(x|\,q)\,,
\qquad\qquad x=\cos\theta\,.
$$
The expression in square brackets factorizes into a product
$({\cal D}_x^{\,q} + 1)({\cal D}_x^{\,q} - 1)$ and the whole equation (7) may be
written as
$$
\left({\cal D}_x^{\,q}\right)^2\, H_n(x|\,q) = q^{-n}\,H_n(x|\,q)\,,   \eqno(8)
$$
where the $q$-difference operator ${\cal D}_x^{\,q}$ is
$$
{\cal D}_x^{\,q} := \,\frac{1}{1- e^{-2{\rm i}\theta}}\,\,e^{\,{\rm i}\,
\ln q^{1/2}\,\partial_{\theta}} + \frac{1}{1- e^{2{\rm i}\theta}}\,
\,e^{-{\rm i}\,\ln q^{1/2}\,\partial_{\theta}}$$
$$\equiv \frac{1}{2{\rm i}\sin\theta}\,\left(e^{{\rm i}\theta}\,
e^{{\rm i}\,\ln q^{1/2}\,\partial_{\theta}} - e^{{-\rm i}\theta}\,
e^{{-\rm i}\,\ln q^{1/2}\,\partial_{\theta}}\,\right)\,.            \eqno(9)
$$

To facilitate ease of clarifying the distinction between ${\cal D}_x^{\,q}$
and the Askey-Wilson divided-difference operator $D_q$, defined by (4), one
may also write (9) in the form
$$
{\cal D}_x^{\,q} \,f(x)= \frac{1-q}{2\sqrt q}\,\frac{1}{\delta_q\,x}\,
\left[e^{-{\rm i}\theta}\,g(q^{1/2}\,e^{{\rm i}\,\theta}) -
e^{{\rm i}\theta}\,g(q^{-1/2}\,e^{{\rm i}\,\theta})\right] $$
$$
= \frac{e^{{\rm i}\theta}g(q^{-1/2}e^{{\rm i}\theta})-e^{-{\rm i}\theta}
g(q^{1/2}e^{{\rm i}\theta})}{e^{{\rm i}\theta}-e^{-{\rm i}\theta}}\,,  \eqno(10)
$$
where $g(e^{{\rm i}\theta})\equiv f(x)$ and $x=\cos\theta$, as before.

For various applications it is important that the $D_q$ and $\widetilde w^{-1}
(x|\,q)\,D_q\,\widetilde w(x|\,q)$ are in fact lowering and raising operators,
respectively, for the continuous $q$-Hermite polynomials $H_n(x|\,q)$ (see, for
example, formulae (3.26.7) and (3.26.9) in \cite{KS}). This circumstance enables
one to interpret a Hilbert space of functions on $[-1,1]$, which are square
integrable with respect to the weight $\widetilde w(x|\,q)$, as a direct sum of
two $su_q(1,1)$-irreducible subspaces $T^+_{1/4}$ and $T^+_{3/4}$, consisting of
even and odd functions, respectively \cite{KD,AS}. So it becomes transparent how
the Askey-Wilson divided-difference operator $D_q$ and the continuous $q$-Hermite
polynomials $H_n(x|\,q)$ are interrelated from the group-theoretic point of view.

An explicit analytic relation between the Askey-Wilson divided-difference operator
$D_q$ and the difference operator ${\cal D}_x^{\,q}$, which surfaces in (8), involves
the so-called {\it averaging difference ope\-rator} ${\cal A}_q$, defined as
$$
\left({\cal A}_q \,f\right)(x)= \,\frac12\,\left( e^{\,{\rm i}\ln q^{1/2}\,
\partial_{\theta}}+ e^{-{\rm i}\ln q^{1/2}\,\partial_{\theta}}\right)f(x)\,
\equiv\,\cos\left(\ln q^{1/2}\,\partial_{\theta}\right)\,f(x)\,.     \eqno(11)
$$
We recall that the averaging operator ${\cal A}_q$ is intimately associated with
the Askey-Wilson operator $D_q$ because the product rule for the latter one is
of the form (see, for example, formula (21.6.4) in \cite{Ism})
$$
D_q\,f(x)\,g(x)\,= \,{\cal A}_q\,f(x)\,D_q\,g(x)\,+\,D_q\,f(x)\,{\cal A}_q\,g(x)\,.
$$
So, from (4), (9) and (11) one concludes that the ${\cal D}_x^{\,q}$ is may be
expressed in terms of the known operators $D_q$ and ${\cal A}_q$ as
$$
{\cal D}_x^{\,q}\,=\,{\cal A}_q \,+\,\frac{1-q}{2{\sqrt q}}\,x\,D_q \,. \eqno(12)
$$

Note that the operator $({\cal D}_x^{\,q})^2$ represents, as equation (8) implies,
an unbounded operator on the Hilbert space $L^2(S^1)$ with the scalar product
$$
\langle g_1,g_2 \rangle =\frac1{2\pi} \int^1_{-1}
\,g_1(x)\,\overline{g_2(x)}\, \widetilde w(x|\,q)\, dx \,,  \eqno(13)
$$
where the weight function $\widetilde w(x|\,q)$ is defined by (2). In view
of (1) the polynomials $p_n(x):=(q^{n+1};q)_\infty^{-1/2} H_n(x|q)$, $n=0,1,2,
\cdots$, constitute an orthonormal basis in this space such that
$\left({\cal D}_x^{\,q}\right)^2 p_n(x)=q^{-n}p_n(x)$. In particular, the
operator $({\cal D}_x^{\,q})^2$ is defined on the linear span ${\cal H}$
of the basis functions $p_n(x)$, which is everywhere dense in
$L^2(S^1)$. We close $\left({\cal D}_x^{\,q}\right)^2$ with respect to the scalar
product (13). Since $\left({\cal D}_x^{\,q}\right)^2$ is diagonal with
respect to the orthonormal basis $p_n(x)$, $n=0,1,2,\cdots$, its closure
$\overline{\left({\cal D}_x^{\,q}\right)^2}$ is a self-adjoint operator
which coincides on ${\cal H}$ with $\left({\cal D}_x^{\,q}\right)^2$. According to
the theory of self-adjoint operators (see \cite{AG}, Chapter 6), we can take
a square root of the operator $\overline{\left({\cal D}_x^{\,q}\right)^2}$.
This square root is a self-adjoint operator too and has the same eigenfunctions
as the operator $\overline{\left({\cal D}_x^{\,q}\right)^2}$ does. We denote this
operator by $\overline{{\cal D}_x^{\,q}}$. It is evident that on the subspace
${\cal H}$ the operator $\overline{{\cal D}_x^{\,q}}$ coincides with the
${\cal D}_x^{\,q}$. That is, the ${\cal D}_x^{\,q}$ is a well-defined operator
on the Hilbert space $L^2(S_1)$ with an everywhere dense subspace of definition.
Moreover, according to the definition of a function of a self-adjoint operator
(see \cite{AG}, Chapter 6), we have $\overline{{\cal D}_x^{\,q}}p_n(x)=
q^{-n/2}p_n(x)$, that is
$$
{\cal D}_x^{\,q}\, H_n(x|\,q)\,\equiv \,\left[{\cal A}_q \,+\,\frac{1-q}
{2{\sqrt q}}\,x\,D_q \,\right]\, H_n(x|\,q)\,= q^{-n/2}\,H_n(x|\,q)\,.   \eqno(14)
$$
Thus, the continuous $q$-Hermite polynomials are in fact governed by a simpler
$q$-difference equation (14) which is, in essence, a factorized form of (8).

This is a place to point out that the first explicit statement of equation (14),
that we know, is in \cite{AAV-B} and \cite{AA}: in the paper \cite{AAV-B}, it
was stated  without proof, whereas in \cite{AA} it was proved by employing the
Rogers generating function
$$
\sum_{n=0}^{\infty}\,\frac{t^n}{(q;q)_n}\,H_n(x|\,q)=
\left(t\,e^{{\rm i}\theta},t\,e^{-{\rm i}\theta};q\right)
_{\infty}^{-1}\,                                        \eqno(15)
$$
for the continuous $q$-Hermite polynomials $H_n(x;q)$ (see \cite{GR},
p. 26) as follows. Apply the $q$-difference operator ${\cal D}_x^{\,q}$
to both sides of generating function (15) to derive that
$$
\sum_{n=0}^{\infty}\,\frac{t^n}{(q;q)_n}\,{\cal D}_x^{\,q}\,
H_n(x|\,q)= {\cal D}_x^{\,q}\,\left(t\,e^{{\rm i}\theta},t\,
e^{-{\rm i}\theta};q\right)_{\infty} ^{-1}$$
$$
= \left(\,q^{-1/2}\,t\,e^{{\rm i}\theta},\,q^{-1/2}\,t\,
e^{-{\rm i}\theta};q\right)_ {\infty}^{-1}\,=
\sum_{n=0}^{\infty}\,\frac{t^n}{(q;q)_n}\,q^{-n/2}\,H_n(x|\,q)\,.
$$
Then equate coefficients of the same powers of $t$ on the extremal
sides above, to obtain the proof that equation (14) is consistent
with the generating function (15). However, it should be noted that
neither \cite{AAV-B} nor \cite{AA} does contain any discussion of
connection between $q$-difference equations (14) and (3) or (7).

In the limit as $q\to 1$ the continuous $q$-Hermite polynomials $H_n(x|\,q)$
are known to reduce to the ordinary Hermite polynomials $H_n(x)$ (see, for
example, \cite{KS}, p.144),
$$
\lim_{q\to 1}\,\kappa^{-n}\, H_n(\kappa\,x|\,q)= H_n(x)\,, \qquad
\quad \kappa := \sqrt{\frac {1-q}{2}}\,.                        \eqno(16)
$$
Hence, if one rescales $x\to\kappa x$ and then lets $q \to1$ in $q$-difference
equations (8) and (14), both of these equations reduce to the same second-order
differential equation
$$
\left( \partial_x^2 - 2x\,\partial_x + 2n\right)\,H_n(x)= 0\,,\qquad
\quad \partial_x\equiv \frac{d}{dx}\,,
$$
for the ordinary Hermite polynomials $H_n(x)$. This fact is an immediate
consequence of the limit property
$$
\lim_{q\to1}\left[\frac{1}{1-q}\left({\cal D}_{\kappa x}^{\,q}\,-\,I\right)
\right]= \frac{1}{2}\left( x - \frac{1}{2}\,\frac{d}{dx}\right)\frac{d}{dx} \eqno (17)
$$
of the $q$-difference operator ${\cal D}_x^{\,q}$, which is straightforward
checked by employing its definition (9) or (10). Note that the rescaling
parameter $\kappa$ in (17) is the same as in (16), whereas $I$ is the identity
operator.

Observe also that by combining (14) and (6) one arrives at the $q$-difference
equation
$$
{\cal D}_x^{1/q}\, H_n(x|\,q)\,{\widetilde w}(x|q) =
q^{-(n+1)/2}\,H_n(x|\,q)\, {\widetilde w}(x|q)\,,                \eqno(18)
$$
which can be viewed as a factorized form of the conventional $q$-difference
equation (3).

In the foregoing exposition up to the present point it has been implied that
$0<q<1$. Of course, the case of $q>1$ can be treated in a similar way. We
briefly state below some explicit formulas for $q>1$ case too. As was noticed
by Askey \cite{A}, one should deal with the case of the continuous $q$-Hermite
polynomials $H_n(x|\,q)$ of Rogers when $q>1$ by introducing a family of
polynomials
$$
h_n(x|\,q):= {\rm i}^{-n}\,H_n({\rm i}\,x|\,q)\,,                \eqno(19)
$$
which are called the continuous $q^{-1}$-Hermite polynomials \cite{IM}.
So the transformation $q \to q^{-1}$ and the change of variables $\theta
= \pi/2-{\rm i}\varphi$ in the $q$-difference equation (14) converts it,
on account of the definition (19), into equation
$$
{\widetilde {\cal D}}_x^{\,q}\, h_n(x|\,q) = q^{n/2}\,h_n(x|\,q)\,,
\quad\quad x=\sinh\varphi\,,                                     \eqno(20)
$$
where the $q$-difference operator  ${\widetilde {\cal D}}_x^{\,q}$ is of
the form
$$
{\widetilde {\cal D}}_x^{\,q} := \frac{1}{2\cosh \varphi}\,\left(e^{\,\varphi}
\,e^{\ln q^{1/2}\,\partial_{\varphi}} + e^{-\varphi}\,e^{-\ln q^{1/2}\,
\partial_{\varphi}}\right)\,.                                         \eqno(21)
$$
One may verify that this $q$-difference equation (20) is in agreement
with the generating function
$$
\sum_{n=0}^{\infty}\,\frac{q^{\,n(n-1)/2}}{(q;q)_n}\,t^{\,n}\,
h_n(\sinh\varphi|\,q) = \left( t\,e^{-\varphi}\,,-t\,e^{\,\varphi};
\,q\right)_{\infty}
$$
for the continuous $q^{-1}$-Hermite polynomials $h_n(x|\,q)$ \cite{IM}.
The proof of (20) follows the same lines as the proof of (14) using the
results of Appendix of the paper \cite{MPL}. Namely, an analogue of the
$q$-difference equation (3) is obtained from the raising and lowering
operators (A.1) and (A.2) in \cite{MPL}. The formula (A.10) in \cite{MPL}
gives an analogue of the relation (6). An analogue of the Hilbert space
$L^2(S^1)$ is constructed in the following way. We take the Hilbert space
$L^2(\mathbb{R})$ with the scalar product determined by the formula
(3.2) in \cite{A}. The polynomials $h_n(x)$, $n=0,1,2,\cdots$, are orthogonal
in this Hilbert space. However, this set of polynomials does not constitute
a basis of $L^2(\mathbb{R})$ since the orthogonality measure in
(3.2) of \cite{A} is not extremal for the continuous $q^{-1}$-Hermite
polynomials. For this reason, we create the closed subspace ${\cal L}$ of
$L^2(\mathbb{R})$ spanned by the polynomials $h_n(x)$, $n=0,1,2,\cdots$.
It is shown in the same way as above that
$({\widetilde {\cal D}}_x^{\,q})^2$ is a bounded self-adjoint operator
on the Hilbert space ${\cal L}$, diagonalizable by the polynomials
$h_n(x|\,q)$, $n=0,1,2,\cdots$. Therefore, the relation (19) holds
for the operator ${\widetilde {\cal D}}_x^{\,q}$.

In conclusion, this short paper should be considered as an attempt to call
attention to a curious fact that the conventional $q$-difference
equation (7) for the continuous $q$-Hermite polynomials $H_n(x|\,q)$ of Rogers
admits factorization of the form $\left[\left({\cal D}_x^{\,q}\right)^2\, - 1
\right]\,H_n(x|\,q)=(q^{-n}- 1)\,H_n(x|\,q)$, where ${\cal D}_x^{\,q}$ is
defined by (9). This circumstance seems to have escaped the notice of all
those with whom we share interests in $q$-special functions.

Finally, since the continuous $q$-Hermite polynomials $H_n(x|\,q)$
occupy the lowest level in the hierarchy of ${}_4\phi_3$ polynomials
with positive orthogonality measures, it is of interest to find out
whether there are instances from higher levels in the Askey $q$-scheme
\cite{KS}, which also admit factorization of an appropriate $q$-difference
equation.

It is well-known that the $q$-different equation (3) is related to the
Hamiltonian for the $q$-oscillator \cite{AS}. So it would be of interest
to look for some insight into the equations (14), (18) and (20) physically.

We are grateful to N. M. Atakishiyev for suggesting to us the problem
and a helpful discussion.

\end{document}